# Dynamics of balanced parentheses, lexicographic series and Dyck polynomials

Gennady Eremin

ergenns@gmail.com

September 18, 2019

**Abstract.** The article deals with a lexicographic order in various sequences. Consider the axiomatic of lexicographic series, based on the properties of the natural numbers. Elements of the set are ordered first the code length; further in each sign range, sorting is performed according to the given order on the alphabet. The sequence of the Dyck words, Dyck series, is analyzed as an example of such lexicographical series. The basis of this series is the dynamics of the Dyck words. We solve the direct and inverse problem of identification of elements of the Dyck series. The polynomial equation on the Dyck triangle is investigated. We give a recursive equation for Dyck polynomials. A matrix of polynomial coefficients is constructed to solve some problems. In conclusion, the reader is offered a software service for identification of the Dyck words with index up to $10^{10}$.

## 1   Introduction

**1.1. Dyck words.** In discrete mathematics, the balanced parentheses, *Dyck words*, are sufficiently known and play an important role [St11]. We present a Dyck word as a system of interrelated elements from an alphabet with two characters. Usually this is a string of left (open) and right (closed) parentheses that are balanced. The system of related parentheses is characterized by the dynamics of Dyck words, *Dyck dynamics*.

For the Dyck word, first, the number of left and right parentheses is the same (the *first rule* of dynamics). And secondly, in every initial subword, the number of right parentheses never exceeds the number of left ones (the *second rule*). Any open parenthesis has a matching closed one and they must be correctly nested.

For any Dyck word of length 2*n*, *Dyck 2n-word* (there are *n* left parentheses and *n* right ones), the second rule of dynamics is equivalent to the known condition for the position $r_i$ of the *i*-th right parenthesis [Ka09]: $2i \leq r_i \leq n + i$, $1 \leq i \leq n$. A group of consecutive Dyck words is also a Dyck word. There are no restrictions to the length of the bracket sequences, so we can talk about the infinity of the set of the Dyck words.

**Example 1.1.** Let's show some variants of balanced brackets and unbalanced ones.
Balanced brackets:  (()),  ()(),  ((()())).
Unbalanced brackets:  )(()),  ()((,  ((()(()).

The number of Dyck 2*n*-words is the *n*th Catalan number (see OEIS A000108):

(1.1) $$c_n = \binom{2n}{n}/(n+1) = \frac{(2n)!}{n!(n+1)!}, \quad n \geq 0.$$

The first Catalan numbers for *n* = 0, 1, 2, … are 1, 1, 2, 5, 14, 42, 132, 429, …

**1.2. Series in mathematics.** Brockhaus and Efron give this definition of a series: "*A series is a sequence of elements composed according to some law. For a given series, there is a law by which you can make as many elements of this series*". According to Brockhaus and Efron, the elements of the series can be numbers, functions and actions. Let's expand this list to include balanced brackets.

Consider as an example the natural numbers, the most famous, the most cited and the simplest numerical sequence:

$$[0,] 1, 2, 3, 4, 5, 6, 7, 8, 9, 10, 11 \ldots$$

A natural series can start with 0 or 1; zero can be a natural number or not. It is logical to call 0 a *runaway zero*. The law (or *consecution function*) by which the series is constructed is very simple: the next element is 1 more than the current one.

There is still no agreement on zero among mathematicians. Some definitions begin the natural numbers with 0 (non-negative integers), whereas others start with 1 (positive integers). Let's try to clarify the situation with this zero; I will give my possibly controversial views on this matter.

In almost all infinite sequences, the element index is a natural number. And either we use a runaway zero or not, depending on there is an element with the index 0 in the sequence or not. We can say that the natural series indexes such sequences. The natural series accompanies sequences, is attached to them.

Obviously, the natural series indexes itself (the *self-indexing* series): every element is its own index. And if we include 0 in the natural series, then this element should automatically be assigned a zero index. Let's write some sequence with direct indexing of elements; any element has a unique (non-repeating) index:

$$[a_0,] a_1, a_2, a_3, a_4, a_5, a_6, a_7, a_8, a_9, a_{10}, a_{11} \ldots$$

An element with index 0 can be or not (again there is a runaway element). In this regard, let us formulate a logical statement.

**Proposition 1.1.** *If we identify a natural series with a family of element indices of some sequence (which is logical) and if in this sequence the initial element has a zero index, then we should consider 0 as a natural number.*

For many sequences, elements are indexed from 0. In numeric series, the initial element can be zero or not. For example, in the Fibonacci sequence, the initial element is 0. For Catalan numbers and Motzkin numbers the zero elements is 1. It is necessary to note such paradox: *many mathematicians admit elements with zero indices in sequences, they explain and justify the presence of such* (*often virtual*) *elements, but at the same time they do not consider 0 a natural number*.

## 2    Lexicographic series

In non-numeric sequences, the lexicographic order (ranking of elements) is always welcome. This order involves strict rules for the construction of elements (objects), as well as some structural principles of the series. Let's build a sequence of non-numeric objects, which will be observed all the characteristic features of the natural numbers. We will call such a sequence a *lexicographic series* or a *lex-seri*es.

**2.1. Axiomatic of lex-series.** What are the formal features of the natural numbers? As you know, integers are not repeated in the natural series. In this regard, we formulate the first axiom of a lex-series.

**Axiom 2.1.** *In a lexicographic series, all objects are unique.*

Some sequences do not satisfy Axiom 2.1, for example, $c_0 = c_1 = 1$.
    Next the natural numbers are arranged (sorted) in order of increasing the code length. The list begins with single-digit integers, followed by double-digit integers, etc. Of course, we can add leading zeros to any integer. This will not change the integer value, but will break the sorting. But this is not usually done, and there is no need. Let's formulate the following axiom of a lex-series.

**Axiom 2.2.** *In a lexicographic series, objects are sorted in ascending order of code length.*

According to Axiom 2.2, lex-series objects are divided into *ranges*, which contain elements with the same code length. Thus, a lex-series is cut into such ranges.
    The natural numbers within each range are also sorted. In the decimal notation, each digit of the alphabet has its own unique weight. The alphabet of the natural numbers is totally ordered as follows:

(2.1) $$0 < 1 < 2 < 3 < 4 < 5 < 6 < 7 < 8 < 9.$$

Obviously, as the weight of a natural number $n$, $wt(n)$, we can consider its index (that is, its value). For example, $wt(0) = 0$, $wt(1) = 1$, $wt(99) = 99$, and so on.
    In accordance with (2.1), integers are sorted within each range. By the way, for Dyck words we have a simple alphabet $A_{dyck} = \{$ '(', ')' $\}$; there are only two possible ordering options: '(' < ')' or '(' > ')'. Both orders can be found in the literature (for example, see [Sa18] and here). Let's formalize the corresponding following axiom of a lex-series.

**Axiom 2.3.** *In each range of a lexicographic series, objects are sorted according to the given totally order on the alphabet.*

The sorting algorithm within the range is simple. When analyzing two elements of the same length, the code signs are compared sequentially in pairs, moving from left to right. And the first mismatch in any pair (this is inevitable, since all objects are unique) allows you to correlate the weights of the compared elements.

**Example 2.4.** Let the totally order '<' be given on the alphabet $A$ of the lex-series $L$ and let $a = a_1a_2…a_k$ and $b = b_1b_2…b_k$ be elements from the *k-range* of $L$. We have $a < b$ if and only if $a_i < b_i$ for the first $i$ where $a_i \neq b_i$. □

So, in a lex-series we have a double sorting of objects. External sorting orders elements along the code length, distributes them to ranges. Internal sorting arranges the elements within the range. Any range has minimum and maximum elements. For example, in the natural series, the minimum element of the $k$-range is 1 and $k–1$ zeros (10…0), the maximum is $k$ nines (99…9). In a lex-series, objects are placed in accordance with their weights; it is convenient and logical to consider the object index as its weight.

Let's analyze the smallest digit 0. Zero is not written at the beginning of integers. The exception is the number 0 itself, if zero is included in the natural series. We are ready to formulate the last axiom of a lex-series about prohibition of placement and replication of the symbol with the minimum weight at the beginning of the code, but one circumstance prevents.

The digits into integers are *free* and practically unrestricted in use. Any sign (except zero) can be repeated at the beginning of the code. At the same time, in the alphabet of Dyck words (here only two signs), both left and right parentheses are not free, but *connected*. In the Dyck $2n$-word, signs can be divided into $n$ pairs; each parenthesis corresponds to some opposite parenthesis.

In the alphabet of the natural numbers there are no associated characters, and in the alphabet of Dyck words there are no free characters. But there is an alphabet $A_{motzkin} = \{\ 0,\ '(',\ ')'\ \}$ of Motzkin words, in which there are connected parentheses (similar to Dyck words) and free 0 (an analogue of zero in the natural numbers). In the Motzkin word, we can put zero anywhere in the code, and even repeat many times at the beginning of the word. Usually zero has minimal weight in the Motzkin word alphabet, and this is logical. In this regard, we formulate the last axiom of a lex-series.

**Axiom 2.5.** *In a lexicographic series, an element with two or more signs in the code cannot begin with a free character that has the least weight in the alphabet.*

Obviously, it is impossible to build the lex-series from Motzkin words. But if we remove leading zeros in Motzkin words and then get rid of duplicates, then we can build the lex-series from such *truncated words*.

**Note 2.6.** The Dyck 2*n*-word is often written as a numeric string of 2*n* binary digits (for example, see [KA09]). In OEIS A063171, the left parenthesis is replaced by 1, the right one is 0. In another sequence (see OEIS A014486), these binary numbers are represented in decimal notation. Seems this way of encoding the Dyck word is illogical.

The sequence of the Dyck words and the sequence of the Motzkin words are very close, since any Dyck word is also a Motzkin word. There are tasks, which included both the Dyck words and Motzkin words. In the alphabet of Motzkin word we have zero; in digitizing the codes we have to keep the sign 0 in the numeric string. In the Motzkin word, there is no need to encode the sign 0 with another digit.

In this regard an interesting way of digitizing, when the Dyck 2*n*-word is a string of *n* 1's and *n* 2's (for example, see Glenn Tesler).

**2.2. Dyck series.** For Dyck words of equal length, let's establish a lexicographical order based on the inequality '(' < ')'. The number of Dyck 2*k*-words is equal to $c_k$. For example, Dyck 6-words form a chain of five elements:

$$((())),\ (()()),\ (())(),\ ()(()),\ ()()().$$ 

In such lexicographic series, or *Dyck series*, 2*k*-range (*k* pairs of parentheses) starts with the single-block word '((…())…)' and ends on the word '()()…()', in which *k* blocks. Let's show the beginning of the Dyck series

(2.2) $\qquad$ (), (()), ()(), ((())), …, ()()(), ((())), …

First we specify a single Dyck word with one pair of parentheses (2-range), followed by two elements with two pairs of parentheses (4-range), etc. In the sequence OEIS A000108, the initial Catalan number is 1 and has an index of 0. And the question arises: what should we do with an empty word of length 0? Whether to include an additional 0-range with a single element (empty word) at the beginning of the Dyck series?

There is no empty word in our Dyck series, and this is almost always the case when listing balanced parentheses (see [La03], page 35). Moreover, you can find a shortened sequence of Catalan numbers without an initial element, for example, as here. And this is understandable; we are usually interested in real bracket sets containing at least one pair of parentheses. In the future, if possible, we will do without an empty word and, accordingly, without the Catalan number with index 0.

For Dyck series, we need to solve the following classical problems:
- For a given Dyck word find the next (preceding) element of the series.
- Generation in lexicographical order of all Dyck words of a given length.
- For a given Dyck word calculate its index in the lex-series.
- Reconstruction Dyck word by its index in the lex-series.

The first two tasks are not difficult and well described in the literature. We refer to the website [Sa18], where there are appropriate algorithms and program codes. In this article, we look at the last two listed tasks.

Recall, all elements of the Dyck series are unique, so each Dyck word occupies a fixed place, i.e. has a specific index. Each natural number corresponds to a single Dyck word. In other words, there is a one-to-one correspondence between Dyck words and the natural numbers. In order to implement this correspondence, we just need to solve two identification problems from the above list.

The index of an arbitrary Dyck word can be given in two ways: (a) an *absolute index* in the Dyck series or (b) a *relative index* – an ordinal number in a given series range. The relation between the absolute index $I$ and the relative index $I_{rel}$ of a given Dyck word of semilength $n$ is obvious:

(2.3)     $I_{abs} = I_{rel} + c_1 + c_2 + \ldots + c_{n-1}, \ I_n \in [1, c_n]$.

Recall that the Dyck series begins with the Dyck word '( )'.

Indexing of the Dyck word and the inverse problem (restoration of the bracket set by a given index) are solved with the involvement of Dyck dynamics. In addition, we refer to this resource where the correct algorithms are given without explanation and without theoretical calculations, in our opinion.

## 3  Dyck dynamics

**3.1. Dyck paths.**  Dyck words are often associated with *Dyck paths* on an integer lattice in the positive quadrant. A Dyck word corresponds to a continuous polyline composed of diagonal steps up (1, 1) and down (1, –1), issuing from the origin and ending on the *x*-axis without crossing the *x*-axis. The left parenthesis is an up-step; the right one is a down-step. The figure shows the path for Dyck word ((()(()()))).

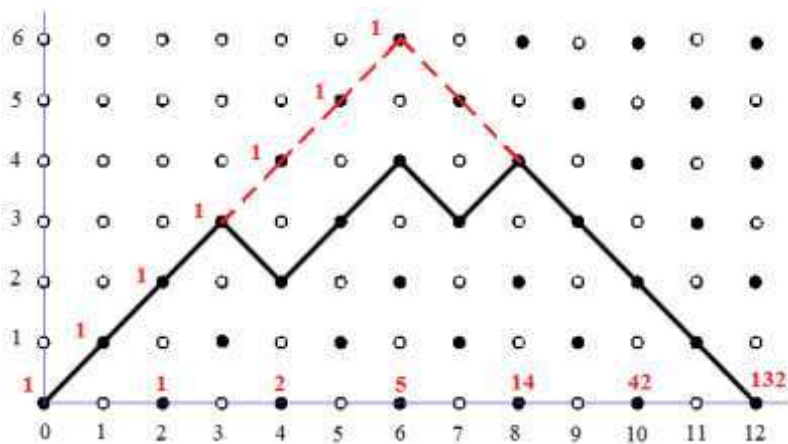

The *x*-axis is the *position* of the parentheses, and *y*-axis is *unbalance* of the parentheses (the number of left parentheses over right ones).  Each left parenthesis increases the height of the polyline, right parentheses reduce the unbalance. The polyline reaches a height of 4,  the maximum height of 6 corresponds to the initial Dyck word ((((((())))))). A triangle with the constructed top (red dotted line) is called a *supporting triangle*.  Each polyline of length 12 does not extend beyond the supporting triangle.

Grid nodes marked with transparent points are unreachable (forbidden) for broken lines. For reachable node, the sum of the coordinates is even. Obviously, the separation of reachable and forbidden nodes is preserved when the supporting triangle is expanded. The number of Dyck paths of length $2n$, *Dyck 2n-path*, is $c_n$; this is the number of Dyck $2n$-paths that can be drawn from $(0, 0)$ to $(2n, 0)$ in the supporting triangle of height $n$. In this regard, node $(12, 0)$ has a red label $c_6 = 132$. On the *x*-axis, five other nodes are similarly marked.

Let's mark the other nodes of the supporting triangle. The label of the reachable node is equal to the number of paths from the origin to this node (label 0 corresponds to each forbidden node). A single ray of ascending links leads to each node of the central diagonal. We marked the diagonal nodes (including the origin) with 1. Specified marks are obvious. With increasing the number of parentheses we get the Catalan numbers on the *x*-axis and we have a chain of 1's diagonally. To understand the internal nodes of the supporting triangle, it is necessary to analyze the mutual links in the broken lines, i.e. the Dyck dynamics.

**3.2. Dynamics equation.** Let's look at a fragment of Dyck path in the next figure. Here we have a down-step that follows the ascending diagonal links. In this case, the last node would have label 1 (as diagonal nodes) if this node was on the *x*-axis. But this is not the case; to the last node from below you can draw an ascending link (shown in green) from the neighboring node with the label, for example, $k$. Then the label of the last node is $1+k$, because all paths ending in both child nodes can be drawn to this node.

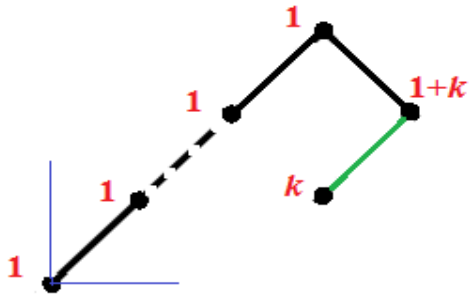

This situation is typical for all internal points: each node can be connect both an ascending and descending link from the child nodes with a smaller abscissa. The label of the inner node is equal to the sum of the labels of the child nodes. Labels below the *x*-axis and above the central diagonal are 0, so the rule is valid for the entire coordinate grid. Thus, for an arbitrary node, its label, *dynamics*, is determined by the following *dynamics equation*:

(3.1) $$d_{i,j} = d_{i-1,j-1} + d_{i-1,j+1}, \quad d_{0,0} = 1, \quad j > 0.$$

Using (3.1) it is easy to calculate the dynamics of all reachable nodes. In the future, we need the *difference dynamics equation*, i.e. calculation of dynamics for the top point, the second summand, in the equation (3.1):

(3.2) $$d_{i,j} = d_{i+1,j-1} - d_{i,j-2}, \quad j > 1.$$

Let's call the equation (3.2) the *dynamics of the top point* or the *top-dynamics*.

**3.3. Dyck triangle.** At any *i*-th position of Dyck word, the unbalance cannot exceed *i*, so the set of reachable nodes forms a triangular matrix, which is known as *Dyck triangle* [MO10]. We will distinguish the Dyck triangle (infinite matrix) and the *supporting triangle*, the initial fragment of the matrix, an isosceles triangle of height *n*, where *n* is the number of pairs of parentheses in a given Dyck word.

The figure below shows a fragment of Dyck triangle; we have identified the supporting triangle for $n = 6$. In the triangle, the outer ascending diagonal (labels

| | | | | | | | | | | | | |
|---|---|---|---|---|---|---|---|---|---|---|---|---|
| 8 | | | | | | | 1 | | 9 | | 54 | |
| 7 | | | | | | 1 | | 8 | | 44 | | |
| 6 | | | | | 1 | | 7 | | 35 | | 154 | |
| 5 | | n=6 | | 1 | | 6 | | 27 | | 110 | | |
| 4 | | | 1 | | 5 | | 20 | | 75 | | 275 | |
| 3 | | 1 | | 4 | | 14 | | 48 | | 165 | | |
| 2 | 1 | | 3 | | 9 | | 28 | | 90 | | 297 | |
| 1 | 1 | 2 | | 5 | | 14 | | 42 | | 132 | | |
| 0 | 1 | 1 | 2 | | 5 | | 14 | | 42 | | 132 | |
| | 0 | 1 | 2 | 3 | 4 | 5 | 6 | 7 | 8 | 9 | 10 | 11 | 12 |

1's) and the outer descending diagonal (labels 1, 6, 20, 48, etc.) contain the same number of nodes. For the point $a = (i, j)$, there is the symmetric point $b = (2n–i, j)$ with its own label. Such a pair we call a *mirror* and denote $a \sim b$.

*Self-symmetric nodes* $(n, j)$ are located along the height. In the figure, these are four central nodes labeled 1, 5, 9, and 5. Can be show that the sum of the label squares of the self-symmetric nodes is $c_6$. Check, $1^2 + 5^2 + 9^2 + 5^2 = 132$.

The Dyck $2n$-path starts at the origin and ends at $(2n, 0)$. If you invert this path, starting at $(2n, 0)$, and recalculate the *inverse dynamics* $đ_{i,j}$, then the points $(i, j)$ and $(2n-i, j)$ exchange labels, i.e. $đ_{i,j} = d_{2n-i,j}$ and $đ_{2n-i,j} = d_{i,j}$.

**Proposition 3.1.** *In the supporting triangle of height n, the number of Dyck 2n-paths passing through the point $(i, j)$ is $d_{i,j} \times đ_{i,j}$.*

*Proof.* The label $d_{i,j}$ is the number of paths leading from the origin to $(i, j)$, while the inverse label $đ_{i,j}$ is the number of paths from the endpoint $(2n, 0)$ to $(i, j)$. Hence, the number of Dyck $2n$-paths passing through the point $(i, j)$ is equal to the product of both labels. This concludes the proof. □

Obviously, self-symmetric nodes have the same forward and reverse labels, that is, $d_{n,j} = đ_{n,j}$. In this regard, the following statement is obvious.

**Corollary 3.2.** *In the n-th column of the Dyck triangle, the sum of the label squares is $c_n$.*

# 4 Identification of Dyck words

**4.1. Computing the index of the Dyck word.** We believe that a certain Dyck word is given, and accordingly we know in what range of the series we will look for a relative index. After calculating the relative index we use (2.3) to find the absolute index of the Dyck word. Calculate the relative index can be a direct search of codes from the beginning of the range or reverse search from the range end (the direction of the search is determined by the initial characters of the Dyck word). This is the easiest way if the length of the code is less than about 20 pairs of parentheses. But when doubling the number of parentheses, this method is unrealistic.

The Dyck triangle contains complete information on the relative positions of the elements in the Dyck series. Using Dyck dynamics, you can calculate the relative index in one code scan from left to right. Let's first describe the algorithm for calculating the relative index of Dyck word.

**Algorithm 4.1.** The Dyck word begins with the left parenthesis, so let's take the initial working index = 1 and here's why. When analyzing parentheses, the working index usually grows, but if the first element of the range is specified (all left brackets are grouped at the beginning of the code), the index does not change until the algorithm end. Here we draw two important conclusions: (a) a left parenthesis of the code does not change the working index and (b) final right parentheses also do not affect the working index.

Therefore, any left parenthesis should be ignored by the algorithm, and the code index depends only on the right parentheses followed by at least one left parenthesis. Each encountered right parenthesis (provided that not all left parentheses are searched) means that the working index should be increased, i.e. it is necessary to make a jump on the Dyck series, skipping some elements with smaller indices. It remains to determine the magnitude of the jump.

The algorithm ends after viewing all the left parentheses of the code. Therefore, to determine the code index, it is enough to fix all the right parentheses preceding the left ones, calculate and sum the jumps along the Dyck series. □

**Example 4.2.** In section 3.1 we considered the Dyck word

$$((\,(\,)\,(\,(\,)\,(\,)\,)\,)\,)\,).$$

The two right parentheses (highlighted in red) precede the left ones; and these two right parentheses define the relative index of Dyck word. In Dyck triangle, the figure shows the corresponding polyline without the last down-steps (the Dyck path is shown in red). In the supporting triangle we marked the central 6th column (labels 1, 5, 9, and 5). Near the 6th column we are interested in two mirrors (4, 4) ~ (8, 4) and (5, 5) ~ (7, 5) (circled in blue).

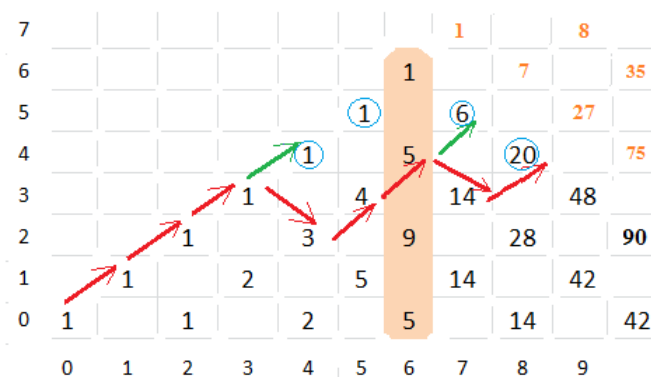

Previously, the working index $I_w$ of the given Dyck word is assumed to be 1. In Dyck path, the first three links are up-steps, and the working index does not change at this area, i.e. $I_w = 1$. But at (3, 3) the path changes direction, and here the working index should be increased, because any path with the 4th up-step (green arrow) from (3, 3) to (4, 4) corresponds to Dyck word with a smaller relative index. The number of such paths is $\breve{d}_{4,4} = d_{8,4} = 20$. There are so many paths with four up-steps at the beginning, and so many Dyck words start with four left parentheses. In short, we need to skip 20 elements of the Dyck series with four left parentheses at the beginning of the code. As a result, it is necessary to perform a jump on the series; the working index becomes equal to $I_w = 1 + 20 = 21$.

The next two links in path are up-steps, so the working index doesn't change. The second break we have at (6, 4), and again it is necessary to perform a jump in the Dyck series. The new green up-step indicates the mirror $(7, 5) \sim (5, 5)$. The analysis shows that it is necessary to skip one more element of the Dyck series, since $\breve{d}_{7,5} = d_{5,5} = 1$. The working index becomes $I_w = 21+1 = 22$.

After passing the last up-step, we come to point (8, 4) on the outer diagonal of the supporting triangle. To fix such output is simple: at any point of such diagonal, the sum of coordinates is $2n$ (in our case 12). Final down-steps do not change the working index; calculation completed.

Thus, the given Dyck word has a relative index of 22 in the 12-range of the series. To obtain an absolute index, you must add the elements of the preceding five ranges. The result is $I_{abs} = 22 + (1+2+5+14+42) = 86$. □

**4.2. Restoration of the Dyck word by its index.** In the previous section, looking through the Dyck word, we analyzed the corresponding path and summarized the labels of selected points in the Dyck triangle in order to obtain the total leap, a relative index in the range of the Dyck series. Now we have the opposite situation: it is necessary to build the Dyck path on the given index. This index should be split into components of the total leap.

Let us be given an absolute index $I_{abs} > 0$. To obtain a relative index we need to subtract the initial Catalan numbers consistently: $I_{abs} - 1 - 2 - 5 - 14 - ...$ and so on until we get a negative number or zero. When subtracting an integer, we skip Dyck words with smaller indices, we move through the series, and we perform a kind of jumps. Obviously, the last positive value is the relative index $I_{rel}$, and the index of the last subtracted Catalan number is the semilength $n$ of the search Dyck word. As a result, $I_{rel} \in [1, c_n]$. If the relative index accepts boundary values, then everything is simple. In case $I_{rel} = 1$, we at once get the initial Dyck word of the $2n$-range '((... ())...)'. If $I_{rel} = c_n$, then we get $n$ consecutive blocks of '( )'.

In General, using the relative index, we need to gradually move from the origin, building step by step the desired Dyck path. At any step, it may be necessary to further reduce the value $I_{rel}$ (additional jumps in the series are performed). This is due to the appearance of right parentheses in the constructed Dyck word. The procedure ends when $I_{rel}$ is zeroed.

**Algorithm 4.3.** We are given a positive index $I_{rel} \leq c_n$. You need to build the suitable Dyck word of the semilength $n$. Any Dyck word starts with the left parenthesis, so we start the lattice path with an up-step from $(0, 0)$ to $(1, 1)$. The first step we entered on the outer diagonal of the supporting triangle; maybe the next steps will follow this diagonal if the given index is small.

In General, we will have to make jumps on the $2n$-range. If the initial word is specified ($I_{rel} = 1$), there will be no jumps. Let's introduce an additional working variable whose value is equal to the required leap, $Jump = I_{rel} - 1$. The value of $Jump$ should be decomposed into components (step-by-step movements over the $2n$-range). Let us ask ourselves: where to move from point $(1, 1)$? Step up to $(2, 2)$ or step down to $(2, 0)$? We get the answer if we compare $Jump$ and the inverse dynamics of point $(2, 2)$. According to the top-dynamics (3.2), the number of Dyck paths that begin with two left parentheses is $\breve{d}_{2,2} = d_{2n-2,2} = c_n - c_{n-1}$.

Three cases are possible:
- $Jump < \breve{d}_{2,2}$. In this case, we cannot perform a corresponding leap to skip all words with two left parentheses at the beginning of the code. So we continue Dyck path by an up-step from $(1, 1)$ to $(2, 2)$ without changing $Jump$. (Perhaps then the value $Jump$ is enough for the leap to skip the words with three parentheses, as $\breve{d}_{3,3} = d_{2n-3,3} = c_n - 2c_{n-1}$.);
- $Jump > \breve{d}_{2,2}$. Now the leap is possible and it must be made. We make a down-step down from $(1, 1)$ to $(2, 0)$, moving to the inner diagonal. At the same time reduce the working variable: $Jump = Jump - \breve{d}_{2,2}$;
- $Jump = \breve{d}_{2,2}$. In this case, the resource $Jump$ is exhausted, all necessary jumps are performed. We have to put the remaining up-steps and then down-steps for Dyck path.

Thus, the procedure is repeated step by step until $Jump > 0$. □

**Example 4.4.** Let's build the Dyck word with $I_{abs} = 1329$. Calculate the relative index
$$I_{rel} = 1329 - 1 - 2 - 5 - 14 - 42 - 132 - 429 = 704.$$
The last thing we subtracted $c_7$, so we will look for the Dyck word in the 16-range. The initial leap is $Jump = 704 - 1 = 703$. Below is a fragment of the supporting triangle; for

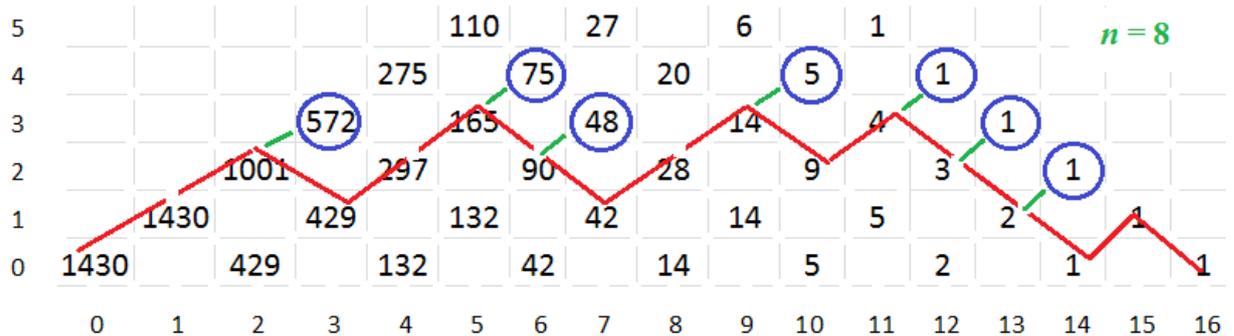

convenience, *reverse labels* (inverse dynamics) are placed in the points. The Dyck

16-path is drawn in red, and the inverse labels of the points that affected the shape of the path are circled in blue.

The first fracture occurred in (2, 2), we made a leap to 572 in the 16-range. As a result, the working value of the leap decreased $Jump = 703 - 572 = 131$. The second fracture was formed in (5, 3), and the next two down-steps give a total shift of $75 + 48 = 123$. Now the working value of the leap equal $Jump = 132 - 123 = 8$. Next we have a shift of 5 at (9.3) and then three shifts of 1 at (11, 3), (12, 2), (13, 1).

As a result, the resource *Jump* is exhausted (zeroed), the polyline goes to the outer descending diagonal, and the last down-step completes the reconstruction of the Dyck path. If we sum up the inverse labels of the marked points (in the blue ring), we get the initial leap of 703. The required Dyck word has the form (()(())(()()))(). □

## 5  Dyck polynomials

**5.1. Polynomial matrix.** The considered algorithms work with Dyck triangle nodes. Node labels are often huge. Here is the 100$^{th}$ Catalan number (57 decimal digits):

$c_{100}$ = 896 519 947 090 131 496 687 170 070 074 100 632 420 837 521 538 745 909 320.

To work with such numbers, information processing acceleration methods are often associated with task logic and data structures.

Recall that in the Dyck triangle, the *x*-axis is the parenthesis positions, and the *y*-axis is the unbalance of the parentheses. For any triangle node the ordinate cannot exceed the abscissa. The node $(i, j)$, $i + j = 2n$, is placed on the *n*-th *isoline*, a descending diagonal with an upper point $(n, n)$ and a lower point $(2n, 0)$. In the figure, the 7$^{th}$ isoline is highlighted in yellow. The start and end labels of the *n*-th

| | | | | | | | | | | | | | | | |
|---|---|---|---|---|---|---|---|---|---|---|---|---|---|---|---|
| 8 | | | | | | | | 1 | | 9 | | 54 | | 273 | |
| 7 | | | | | | | 1 | | 8 | | 44 | | 208 | | 910 |
| 6 | | | | | | 1 | | 7 | | 35 | | 154 | | 637 | |
| 5 | | | | | 1 | | 6 | | 27 | | 110 | | 429 | | 1638 |
| 4 | | | | 1 | | 5 | | 20 | | 75 | | 275 | | 1001 | |
| 3 | | | 1 | | 4 | | 14 | | 48 | | 165 | | 572 | | 2002 |
| 2 | | 1 | | 3 | | 9 | | 28 | | 90 | | 297 | | 1001 | |
| 1 | 1 | | 2 | | 5 | | 14 | | 42 | | 132 | | 429 | | 1430 |
| 0 | 1 | | 1 | | 2 | | 5 | | 14 | | 42 | | 132 | | 429 |
| | 0 | 1 | 2 | 3 | 4 | 5 | 6 | 7 | 8 | 9 | 10 | 11 | 12 | 13 | 14 | 15 |

isoline are obvious:

$d_{n,n} = 1$, $d_{n+1,n-1} = n$, $d_{2n,0} = d_{2n-1,1} = c_n$, $d_{2n-2,2} = c_n - c_{n-1}$, $d_{2n-3,3} = c_n - 2c_{n-1}$,
$d_{2n-4,4} = c_n - 2c_{n-1} - (c_{n-1} - c_{n-2}) = c_n - 3c_{n-1} + c_{n-2}$.

The points of the $n$-th diagonal can be calculated on the basis of the difference dynamics equation, let's repeat it:

(5.1) $$d_{i,j} = d_{i+1,j-1} - d_{i,j-2}, \quad j > 1.$$

Let's use the (5.1) and write the initial lower points of the $n$-th diagonal:

$d_{i,0} = c_n$, $n = i/2$ ($i$ is even only);
$d_{i,1} = c_n$, $n = (i+1)/2$ ($i$ is odd only);
$d_{i,2} = c_n - c_{n-1}$, $n = (i+2)/2$ ($i$ is even again);
$d_{i,3} = c_n - 2c_{n-1}$, $n = (i+3)/2$ ($i$ is odd again);
$d_{i,4} = c_n - 3c_{n-1} + c_{n-2}$, $n = (i+4)/2$ ($i$ is even again);
$d_{i,5} = c_n - 4c_{n-1} + 3c_{n-2}$, $n = (i+5)/2$ ($i$ is odd again).

For the other points we need to look for rules for constructing such polynomials.

In the given equalities, on the right we sum up the Catalan numbers, which are indexed by the variable $n$. The coefficients for the Catalan numbers are determined by the unbalance $j$ (the line number in Dyck triangle). In these equalities, the variable $i$ can take certain values (even or odd). Obviously, these equalities can be written as a generalized polynomial

(5.2) $$d_{i,j} = p_j(n), \quad n = (i+j)/2.$$

In the difference equation of dynamics (5.1) the points $d_{i,j}$ and $d_{i+1,j-1}$ are located on the common $n$-th isoline, and $d_{i,j-2}$ is located on the adjacent $(n-1)$-th isoline. Obviously (5.1) and (5.2) are equivalent to the *recursion formula*

(5.3) $$p_j(n) = p_{j-1}(n) - p_{j-2}(n-1), \quad j > 1.$$

The resulting recursion is easy to check, calculate the next 6th polynomial:

$$d_{i,6} = p_6(n) = p_5(n) - p_4(n-1) = (c_n - 4c_{n-1} + 3c_{n-2}) - (c_{n-1} - 3c_{n-2} + c_{n-3})$$
$$= c_n - 5c_{n-1} + 6c_{n-2} - c_{n-3}, \quad n = (i+6)/2, \text{ where } i \text{ is only even.}$$

Let's call the polynomials $p_j(n)$ *Dyck polynomials* and recursion (5.3) a *polynomial equation* or *dynamics of Dyck polynomials*. Note that earlier we understood by *n* a fixed value, the Dyck word semilength. In (5.3) *n* is a parameter of polynomials; this variable is not explicitly represented in the Dyck triangle, $n = (i+j)/2$. Accordingly, the variable $i$ is not explicitly represented in polynomials, $i = 2n - j$.

The coefficients of the Dyck polynomials can be packed into a matrix, the *polynomial matrix*, which takes an almost triangular shape. In the figure below, the $x$-axis shows the unbalance $j$; along the $y$-axis, starting with some conditional $n$, the indices of the Catalan numbers follow in descending order. The $j$-th column contains the coefficients of the $j$-th polynomial. For example, let us repeat the 6$^{th}$ polynomial: $p_6(n) = c_n - 5c_{n-1} + 6c_{n-2} - c_{n-3}$.

|       |   |    |    |    |    |    |    |    |     |     |     |      |      |      |      |      |      |
|-------|---|----|----|----|----|----|----|----|-----|-----|-----|------|------|------|------|------|------|
| n-8   |   |    |    |    |    |    |    |    |     |     |     |      |      |      |      | 1    | 9    |
| n-7   |   |    |    |    |    |    |    |    |     |     |     |      |      |      | -1   | -8   | -36  | -120 |
| n-6   |   |    |    |    |    |    |    |    |     |     |     |      | 1    | 7    | 28   | 84   | 210  | 462  |
| n-5   |   |    |    |    |    |    |    |    |     |     | -1  | -6   | -21  | -56  | -126 | -252 | -462 | -792 |
| n-4   |   |    |    |    |    |    |    |    | 1   | 5   | 15  | 35   | 70   | 126  | 210  | 330  | 495  | 715  |
| n-3   |   |    |    |    |    |    | -1 | -4 | -10 | -20 | -35 | -56  | -84  | -120 | -165 | -220 | -286 | -364 |
| n-2   |   |    |    |    | 1  | 3  | 6  | 10 | 15  | 21  | 28  | 36   | 45   | 55   | 66   | 78   | 91   | 105  |
| n-1   |   | -1 | -2 | -3 | -4 | -5 | -6 | -7 | -8  | -9  | -10 | -11  | -12  | -13  | -14  | -15  | -16  |
| n     | 1 | 1  | 1  | 1  | 1  | 1  | 1  | 1  | 1   | 1   | 1   | 1    | 1    | 1    | 1    | 1    | 1    | 1    |
|       | 0 | 1  | 2  | 3  | 4  | 5  | 6  | 7  | 8   | 9   | 10  | 11   | 12   | 13   | 14   | 15   | 16   | 17   |

Let's build the polynomial matrix. The bottom line is 1's, the remaining lines are initially zeroed. Then the coefficients are calculated line by line according to the formula:

(5.4) $$a_{j,n-k} = a_{j-1,n-k} - a_{j-2,n-k+1}, \quad j > 1, k > 0.$$

Equality (5.4) follows directly from (5.3). The construction of the polynomial matrix begins with the point $a_{2,n-1} = a_{1,n-1} - a_{0,n} = 0 - 1 = -1$.

It is easy to see, on the outer quasi-diagonal, the initial elements of the lines (the upper points of even columns) are an alternating series of units:

$$a_{2k,n-k} = (-1)^k, \quad k > 0.$$

For explanations, the figure highlights three points associated with equality (5.4): $a_{12,n-3} = a_{11,n-3} - a_{10,n-2} = -56 - 28 = -84$.

**Example 5.1.** Let's calculate the dynamics of some nodes of Dyck triangle. The coefficients for the polynomials will be chosen from the polynomial matrix (and initial Catalan numbers can be found in [OEIS A000108](#)).

1) Check the node label (15, 7). Recall, the dynamics of $d_{15,7}$ indicates the number of initial fragments of the Dyck words with unbalance of 7 in position 15. Point (15, 7) is on the isoline $n = (15+7)/2 = 11$. From the coefficients of the 7$^{th}$ column of the matrix, we make the polynomial $d_{15,7} = p_7(11) = c_{11} - 6c_{10} + 10c_9 - 4c_8 = 58786 - 6\times16796 + 10\times4862 - 4\times1430 = 910$.

2) Calculate the dynamics for (32, 10). The coordinates give us the 21$^{st}$ isoline; for the calculation we choose the 10$^{th}$ column of the matrix: $d_{32,10} = p_{10}(21) = c_{21} - 9c_{20} + 28c_{19} - 35c_{18} + 15\,c_{17} - c_{16} = 64\,512\,240$.

3) Make a polynomial for (132, 10). The node is on 71$^{st}$ isoline, but we still choose the 10$^{th}$ column: $d_{132,10} = p_{10}(71) = c_{71} - 9c_{70} + 28c_{69} - 35c_{68} + 15\,c_{67} - c_{66}$. For reference:
$c_{66} = 5\,632\,681\,584\,560\,312\,734\,993\,915\,705\,849\,145\,100$,
$c_{67} = \phantom{0}22\,033\,725\,021\,956\,517\,463\,358\,552\,614\,056\,949\,950$,
$c_{68} = \phantom{0}86\,218\,923\,998\,960\,285\,726\,185\,640\,663\,701\,108\,500$,
$c_{69} = 337\,485\,502\,510\,215\,975\,556\,783\,793\,455\,058\,624\,700$,
$c_{70} = 1\,321\,422\,108\,420\,282\,270\,489\,942\,177\,190\,229\,544\,600$,

$c_{71} = 5\ 175\ 569\ 924\ 646\ 105\ 559\ 418\ 940\ 193\ 995\ 065\ 716\ 350$.
Answer: $d_{132,10} = 39\ 575\ 872\ 930\ 789\ 889\ 398\ 293\ 766\ 300\ 107\ 613\ 200$. □

**5.2. Dyck wedge.** The Dyck triangle is bounded above by the diagonal (nodes labeled 1's) and below by the x-axis (the Catalan numbers). The dynamics of an arbitrary node can be calculated using (3.1) moving from the origin. This procedure can be lengthy and laborious if the node is far from the origin. Dyck polynomials allow calculations with fewer operations if the abscissa of a given node is significantly greater than its ordinate; however, in this case we have to deal with more "long" numbers.

In the Dyck triangle, let's fix the point $x = (i_x, j_x)$. We assume that $x$ is inside the triangle, in other words, $i_x > j_x > 0$ (otherwise the dynamics $d_x$ is obvious). The point $x$ is located on the isoline $n_x = (i_x+j_x)/2$; the isoline "cuts off" the isosceles supporting triangle with the boundary points $(0, 0)$, $(n_x, n_x)$, $(2n_x, 0)$ and the base $2n_x$. Consider three groups of points, three segments that are related (or not) to $x$.

Well, the dynamics $d_x$ can be calculated moving from the origin. Nodes that are used in such calculations we will include in the first *main segment*. The point $x$ is also included in the main segment and has a maximum abscissa there.

Dynamics calculation can be performed from $(2n_x, 0)$, rising along the isoline up to $x$ in accordance with (5.1). Nodes that are used in such case we will include in the second *difference segment*. The point $x$ is included in the difference segment and has a maximum ordinate there. *Thus, the main and difference segments have a common node-joint x.*

In the Dyck triangle, many nodes are often not used in the calculation of $d_x$, and these nodes form the third *dead segment*. The third segment is sometimes comparable to the other segments. Let's do a further analysis on a specific example.

**Example 5.2.** The figure below shows the node (15.9) with a dynamics of 350. The supporting triangle (colored nodes) is bounded by $12^{th}$ isoline with extreme nodes (12, 12) and (24, 0). The main segment (blue nodes) is a rotated trapezoid with a base that coincides with the ascending diagonal of the supporting triangle. The difference segment (dark nodes) is a triangle. The side of the main segment and the difference segment diagonal form the 12$^{th}$ isoline.

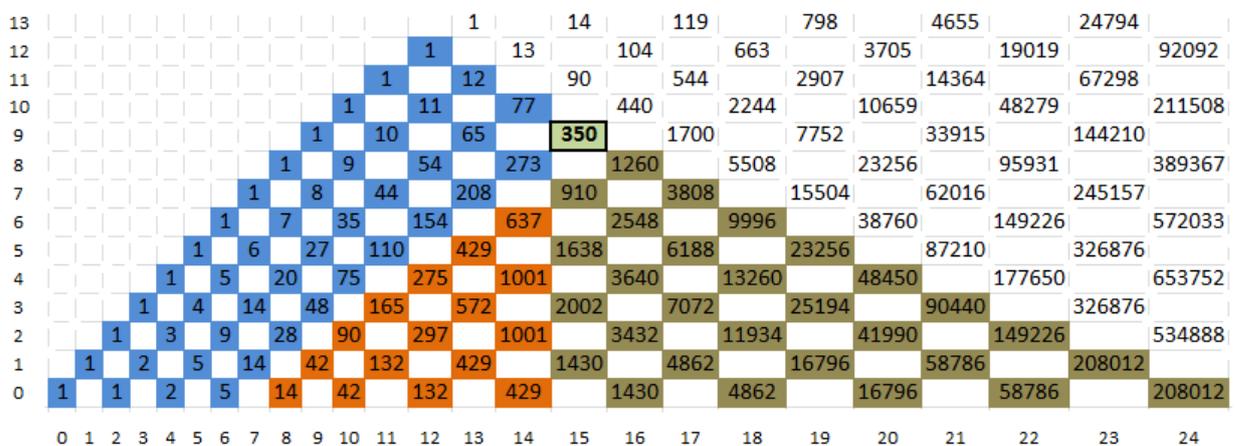

The third segment ("rusty" nodes) is embedded between the main and difference segments. These are the nodes that are not used to calculate $d_{15,9}$; this is the dead zone for (15, 9). The third segment resembles a wedge that corrodes (dissects) the supporting triangle; let's call it *Dyck wedge*. □

In general, for $x = (i_x, j_x)$, the Dyck wedge is placed between the $i_x$-th column and the diagonal that connects $x$ and $(i_x - j_x, 0)$. The wedge top is placed at $(i_x - 1, j_x - 3)$, so the wedge exists if $j_x > 2$. If you rise diagonally from $(2n_x, 0)$, the difference segment gradually expands; the wedge appears starting at $(2n_x - 3, 3)$. At $(n_x, n_x)$ the main segment degenerates into an outer diagonal of the supporting triangle. In this case, the Dyck wedge and the difference segment reach maximum sizes.

# 6 Online software service

The reader is offered a small software service that works online. This section describes two programs for identification of balanced parentheses in the Dyck series. Direct and inverse identification of series elements is solved; the first program calculates the index of the given Dyck word, the second program restores the balanced brackets for the given index.

Data processing is performed using the *modified Dyck triangle*. To do this, in a normal Dyck triangle, unreachable nodes are removed; as a result, we work in a packed triangular array in which the isoline number is plotted along the *x*-axis (see details here). Huge numbers have to be processed, so almost all modules work with "long arithmetic". Programs run on the client side (browser). This means that only HTML, CSS, and JS components are used.

For the given Dyck word, the index in the series is calculated in this program (direct identification problem). First, the correctness of the bracket set is checked, and then the program displays the range number and both indexes – absolute and relative.

Reconstruction of the Dyck word is performed by its absolute index (inverse identification problem). The program displays two indexes (absolute and relative), the range number, and a set of parentheses.

**Acknowledgements.** I would like to thank and express his sincere gratitude to Bruce Sagan (Michigan State University, USA) for helpful comments.

Gzhel State University, Moscow, 140155, Russia

http://www.en.art-gzhel.ru/